\documentclass[MSNbibl,number,citesort,seceqn,dvips]{arxbj}
\usepackage{upgreek}

\aid{0}
\volume{18}
\issue{3}
\pubyear{2012}
\firstpage{823}
\lastpage{835}
\doi{10.3150/11-BEJ363}

\makeatletter
\newcommand{\eqref}[1]{(\ref{#1})}

\DeclareMathAlphabet{\mathbbl} {U}{mt2hrb}{m}{n}
\SetMathAlphabet{\mathbbl}{bold}{U}{mt2hrb}{b}{n}

\newtheorem{proposition}{Proposition}[section]
\newtheorem{theorem}{Theorem}[section]

\newproclaim{definition}{Definition}[section]
\newremark{rem}{Remark}[section]
\newremark{remark}{Remark}
\newremark{example}{Example}

\newcommand{\la}{\lambda}
\newcommand{\sRp}{[0,\infty)}
\newcommand{\sN}{\mathbb{N}}
\newcommand{\sQ}{\mathbb{Q}}
\newcommand{\afrac}[2]{#1/#2}
\newcommand{\kla}[1]{(#1)}
\newcommand{\abs}[1]{|#1|}
\newcommand{\esym}{\mathbbl{E}}
\newcommand{\psym}{{\mathbbl P}}
\newcommand{\PP}[1]{\psym(#1)}
\newcommand{\fexp}{\operatorname{exp}}
\newcommand{\veps}{\varepsilon}
\newcommand{\ff}{L}
\newcommand{\arro}{\rightarrow}
\newcommand{\dto}{\stackrel{d}{\arro}}
\newcommand{\disas}{\stackrel{d}{\sim}}
\newcommand{\gi}[1]{#1^{-1}}
\newcommand{\lap}{\psi}
\newcommand{\imp}{($\Rightarrow$) }
\newcommand{\pmi}{($\Leftarrow$) }
\newcommand{\ii}[3]{\zeta_{#1}(#2,#3)}
\newcommand{\tnu}{\overline{\nu}}
\newcommand{\eat}[1]{#1^\ast}
\newcommand{\yy}{\eat Y}
\newcommand{\xix}{\eat X}
\newcommand{\fix}{\eat F}
\newcommand{\hix}{\eat H}
\newcommand{\pto}[1]{{\mathcal P}_{#1}}
\newcommand{\dick}{\rho}
\makeatother

\begin{document}
\begin{frontmatter}

\title{On the small-time behavior of subordinators}
\runtitle{On the small-time behavior of subordinators}

\begin{aug}
\author[1]{\fnms{Shaul K.} \snm{Bar-Lev}\thanksref{1}\ead[label=e1]{barlev@stat.haifa.ac.il}},
\author[2]{\fnms{Andreas} \snm{L\"opker}\corref{}\thanksref{2}\ead[label=e2]{lopker@hsu-hh.de}}
\and
\author[3]{\fnms{Wolfgang} \snm{Stadje}\thanksref{3}\ead[label=e3]{wolfgang@mathematik.uos.de}}
\runauthor{S.K. Bar-Lev, A. L\"opker and W. Stadje}
\address[1]{Department of Statistics, University of
Haifa, Haifa 31905, Israel.
\printead{e1}}
\address[2]{Department of
Economics and Social Sciences, Helmut Schmidt
University Hamburg, 22043 Hamburg,
Germany.
\printead{e2}}
\address[3]{Department of Mathematics and Computer Science,
University of Osnabr\"uck,
49069 Osnabr\"uck, Germany.
\printead{e3}}
\end{aug}

\received{\smonth{11} \syear{2010}}
\revised{\smonth{2} \syear{2011}}

%
\begin{abstract}
We prove several results on the behavior near $t=0$ of $Y_t^{-t}$ for
certain $(0,\infty)$-valued stochastic processes $(Y_{t})_{t> 0}$. In
particular, we show for L\'evy subordinators that the Pareto law on
$[1,\infty)$ is the only possible weak limit and provide necessary and
sufficient conditions for the convergence. More generally, we also
consider the weak convergence of $tL(Y_t)$ as $t\to0$ for a decreasing
function~$L$ that is slowly varying at zero. Various examples
demonstrating the applicability of the results are presented.
\end{abstract}

%

\end{frontmatter}
%

\section{Introduction}

We consider the behavior near $t=0$ of a stochastic process $(Y_{t})_{t>
0}$ with values in $(0,\infty)$. Let $F_{t}(y)={\mathbbl P}(Y_{t}\leq y)$
and $\psi_{t}(u)={\mathbbl E}(\mathrm{e}^{-uY_{t}})$ be the\vspace*{-1.5pt} distribution function
and the Laplace--Stieltjes transform (LST) of $Y_{t}$ and let
$\stackrel{d}{\rightarrow}$ denote convergence in distribution.

We start with the following observation from \cite{SE}, which is not
difficult to prove. It states that the convergence of $Y_{t}^{-t}$ to some
nondegenerate random variable (r.v.) with distribution function $F^*$
is equivalent to the weak convergence
of the distribution function $u\mapsto1-\psi_{t} (u^{1/t} )$
to $F^{\ast}$.

\begin{proposition}[(see \protect\cite{SE})]\label{be} Assume that
$F_{t}(0)=0$ for all $t$ and let $Y^{\ast}$ be a
r.v. with distribution function $F^{\ast}$
which is not concentrated at one point. Then $Y_{t}^{-t}\stackrel
{d}{\rightarrow}Y^{\ast}$ as $t\rightarrow0$ if
and only if $\psi
_{t} (u^{1/t} )\rightarrow1-F^{\ast}(u)$ as $t\rightarrow0$
at all
continuity points $u$ of $F^{\ast}$.
\end{proposition}

In \cite{SE}, the applicability of Proposition \ref{be} to various
examples was demonstrated. In these examples, the limiting distribution
$F^{\ast}$
turned out to be either a Pareto law with support $[1,\infty)$, or
a mixture of such a Pareto law and a point mass at $%
1$, or an exponential law (possibly shifted to the right).
In general, any distribution on $(0,\infty)$ can occur as~$F^*$ (take
$Y_t = (Y^*)^{-1/t}$),
but it is a challenging question which $F^*$ appear as limits of
`reasonable' processes $Y_t^{-t}$.

In this paper, we study the case when $\psi_{t}(u)=\psi(u)^t$ for
some LST $\psi$; of course this means that $F$ is infinitely divisible,
and the $(0,\infty)$-valued process $(Y_{t})_{t\geq0}$ (with
$Y_{0}\equiv0$) can then be
interpreted as an increasing L\'{e}vy process (a \textit
{subordinator}) with
Laplace exponent $\Phi(u)$ defined by ${\mathbbl
E}(\mathrm{e}^{-uY_{t}})=\mathrm{e}^{-t\Phi(u)}$.

In \cite{SL}, it is proved that for a subclass of exponential
dispersion models (cf. \cite{Jorg}) generated by an infinitely
divisible probability measure $\mu$ on $[0,\infty)$ and associated
with an
unbounded L\'{e}vy measure $\nu$ satisfying $\nu((x,\infty))\sim
-\gamma
\log x$ as $x\to0$, the limit $F^{\ast}$ is a~Pareto type law
supported on $[1,\infty)$.
Our main result below shows that this is indeed the only limit law that
can occur for any subordinator. We
also give several necessary and sufficient conditions for this
convergence to occur.
Combining subordinators and fixed r.v.'s one obtains mixtures of a
Pareto law and the point mass at 1.

The results presented in this paper enable an approximation of the
distribution of $Y_{t}$ for relatively small values of $t$. Note that while
the distribution of $Y_{t}$ can be quite complex, the specific limiting
Pareto law is rather simple to handle. Such numerical approximation
aspects for various distributions $F_{t}$ are subject of future investigations.

This paper is organized as follows. Some preliminary results are presented
in Section \ref{s2}. Under rather mild conditions on the behavior of
$\psi
_{t}$ (which are satisfied for subordinators),
it is shown that ${Y^{\ast}\geq1}$ almost surely. Some other
straightforward results concerning the limiting behavior of products and
sums of stochastic processes are also presented. In Section\vspace*{-1pt} \ref{s4}, we
present necessary and sufficient criteria for the convergence
$Y_{t}^{-t}%
\stackrel{d}{\rightarrow}Y^{\ast}$ for
subordinators in terms of their
characteristics. We also provide an alternative proof of the result of
\cite{SL}.
Section \ref{sec4} presents several applications. In Section \ref{s3}, we consider the
problem under which conditions $%
tL(Y_{t})$ converges weakly as $t\rightarrow0$, if $L$ is a slowly varying
decreasing function with $\lim_{x\rightarrow0}L(x)=\infty$. Clearly,
our original
question concerns the special case $L(x)=-\log x$.

\section{Preliminary results}\label{s2}

Our first result deals with the limiting variable $Y^\ast$. Under a suitable
monotonicity condition on~$\psi_t$, it follows easily that $Y^\ast
\geq1$.

\begin{proposition}\label{pro1}
Suppose that $Y_t^{-t}\stackrel{d}{\rightarrow} Y^\ast$ and that
there are $s>0$ and $y>0$ such that~$\psi_t(u)$
is decreasing in $t$ for $t\in[0,s]$ and $u\in[0,y]$. Then $Y^\ast
\geq1$
almost surely.
\end{proposition}

\begin{pf} Let $u<1$. Then $u^{1/t}$ converges to $0$ as $t\to0$
and if $t>\max\{s,\log u/\log y\}$ then $1\geq\lap_t(u^{1/t})\geq
\lap_s(u^{1/t})$ and $\lap_s(u^{1/t})$ converges to $1$ as $t\to0$.
Since additionally $\lap_t(u^{1/t})\to1-\fix(u)$, it follows that
$\fix(u)=0$ for all $0\leq u<1$.
\end{pf}

\begin{example}[(Stable densities)]
Let $\psi
_{t}(u)=\exp(-au^{t})$, $a>0,$ be the LST of the positive stable
density of type $t\in(0,1)$.
The family $(\psi_t)_{t>0}$ does not satisfy the condition stated in
Proposition~\ref{pro1}. In fact, $\psi_t(u^{1/t}) =1 -\mathrm{e}^{-au}$ for
all $t$ so
that $F^*$ is the exponential distribution with mean
$1/a$, whose support is $[0,\infty)$.
Consider however the distributions belonging to the natural exponential
families generated by
these positive stable densities (with canonical parameter $\theta>0$).
They have LST's $\psi_{t}(u;\theta
)=\exp\{ -a[ (\theta+u)^{t}-\theta^{t}] \}
$ and thus
satisfy the condition of Proposition~\ref{pro1}; in this case it is easily
checked that
$F^*$ has the shifted exponential density $a\mathrm
{e}^{-a(x-1)}1_{(1,\infty
)}(x)$; see also Example 2.5(iii) in \cite{SE}.
\end{example}

\begin{proposition}
\label{pro2} Let $(Y_{i,t})_{t>0}$, $i\in\{1,2\}$, be two independent
families of positive
r.v.'s. If $%
Y_{i,t}^{-t}\stackrel{d}{\rightarrow}Y_{i}^{\ast
}$ as $%
t\rightarrow0$ for $i=1,2$, then
\[
(Y_{1,t}Y_{2,t})^{-t} \stackrel{d}{\rightarrow}%
Y_{1}^{\ast}Y_{2}^{\ast} \quad\mbox{and}\quad(Y_{1,t}+Y_{2,t})^{-t}
\stackrel{d}{\rightarrow} \min\{Y_{1}^{\ast},Y_{2}^{\ast
}\}.
\]
\end{proposition}

\begin{pf} The convergence $(Y_{1,t}Y_{2,t})^{-t}\dto\yy_1 \yy_2$
is trivial.
To prove the second assertion, note that
\[
\lap_{1,t} (u^{1/t} )\lap_{2,t} (u^{1/t} )\to\PP
{\yy_1>u}\PP{\yy_2>u}=\PP{\min\{\yy_1,\yy_2\}>u}
\]
for every $u$ that is a common point of continuity of the functions
$u\mapsto\PP{ Y_i^* >u}$, $i=1,2$.
But $\lap_{1,t}\lap_{2,t}$ is the LST of the sum $Y_{1,t}+Y_{2,t}$
and the result follows immediately from
Proposition~\ref{be}.
\end{pf}

In particular, suppose that $a_{t},b_{t}>0$ are positive
functions with $a_{t}\sim a^{-1/t}$ and $b_{t}\sim b^{-1/t}$ as $%
t\rightarrow0$, with some constants $a,b>0$. Then $Y_{t}^{-t}\stackrel
{d}{\rightarrow}Y^{\ast}$ implies
\[
(a_{t}Y_{t}+b_{t})^{-t}\stackrel{d}{\rightarrow
}\min
\{aY^{\ast},b\}.
\]

%

\section{Small-time behavior of L\'evy subordinators}
\label{s4}

\subsection{The main result}

Let $(Y_{t})_{t\geq0}$ be a subordinator, that is, an increasing L\'{e}
vy process with $Y_0\equiv0$ (see Chapter III in \cite{bertoin}). We
assume that $Y_{t}$ has
no drift, so that the process has the L\'{e}vy--Khintchine-representation $%
\psi_{t}(u)\equiv{\mathbbl E}(\mathrm{e}^{-uY_{t}})=\mathrm
{e}^{-t\Phi(u)}$, where the
Laplace exponent $\Phi$ is given by
\[
\Phi(u)=\int_{0}^{\infty}(1-\mathrm{e}^{-ux})\,\mathrm{d}\nu(x).
\]
Here $\nu$ is the L\'{e}vy measure with support $[0,\infty)$,
satisfying $%
\overline{\nu}(x)\equiv\int_{x}^{\infty}\mathrm{d}\nu(u)<\infty
$ and $%
\int_{0}^{x}u\,\mathrm{d}\nu(u)<\infty$ for all $x>0$. In what
follows, we
write $%
Y=Y_{1} $ and $F(x)={\mathbbl P}(Y\leq x)$.

It is known that a driftless subordinator $Y_{t}$ tends to zero
sub-linearly as $t\rightarrow0$, that is, almost surely, $Y_{t}/t$
tends to
zero as $t\rightarrow0$ (Proposition 8 in \cite{bertoin}). Moreover,
if $h(t)$ is an increasing function such that $h(t)/t$
is also increasing, then (see \cite{bertoin}, Theorem 9)
\[
\mbox{either }\quad\lim_{t\rightarrow0} \frac{Y_{t}}{h(t)}=0\qquad
\mbox{a.s.} \quad\mbox{or}%
\quad\operatorname{lim\,sup}\limits_{t\rightarrow0} \frac
{Y_{t}}{h(t)}=\infty\qquad\mbox{a.s.}
\]
L\'evy processes in general possess the small-time ergodic property
\[
\lim_{t\to0} t^{-1} \esym(f(Y_t)) = \int f(x)\,\mathrm{d}\nu(x)
\]
for bounded continuous functions $f$ vanishing in a neighborhood of the
origin (\cite{Sato}, Corollary~8.9). Letting $P_tf(x)=\esym
(f(Y_t)|Y_0=x)$ and $A$ the infinitesimal generator of the Markov
process $Y_t$, this is nothing else than saying that $P_tf(0)\approx
f(0)+t\cdot Af(0)$ as $t\to0$.

We investigate the limiting behavior of $Y_{t}^{-t}$ as $t\rightarrow0$.
Since $\psi_{t}(u)=\psi(u)^t$ is decreasing in $t$ for fixed $u<1$,
it is an
immediate consequence of Proposition \ref{pro1} that the limit in
distribution, if it exists, will be concentrated on $[1,\infty)$. The
Pareto law $\mathcal{P}_{\gamma}$
with parameter $\gamma>0$ has the distribution function
\[
\Pi_{\gamma}(x)=
(1-x^{-\gamma}) 1_{[ 1,\infty)}(x).
\]

\begin{theorem}\label{ttt}
Let $Y^{\ast}$ be a positive r.v.
which is not concentrated at one point. Let $F^{\ast}(x)={\mathbbl P}%
(Y^{\ast}\leq x)$ be its distribution function. Then the following
statements are
equivalent:

\begin{enumerate}[(S1)]
\item[(S1)] \hypertarget{a1} $Y_t^{-t}\stackrel{d}{\rightarrow}
Y^\ast$ as $t\to0$.
\item[(S2)] \hypertarget{a2} $t\Phi(u^{1/t})\to- \log(1-F^\ast
(u))$ as $t\to0$,
for all continuity points $u$ of $F^\ast$.
\item[(S3)] \hypertarget{a3} $Y_t^{-t}\stackrel{d}{\rightarrow} {%
\mathcal{P}}_{\gamma}$ as $t\to0$ \mbox{ for some } $\gamma>0$.
\end{enumerate}
Furthermore, for any $\gamma>0$ the following statements are equivalent:
\begin{enumerate}[(S4)]
\item[(S4)] \hypertarget{a4} $Y_t^{-t}\stackrel{d}{\rightarrow}
{\mathcal{P}}_{\gamma}$ as $t\to0$.

\item[(S5)] \hypertarget{a5} $\Phi(s)/\log s\to\gamma$ as $s\to
\infty$.

\item[(S6)] \hypertarget{a6} $\log F(x)/\log x\to\gamma$ as $x\to0$.

\item[(S7)] \hypertarget{a7} $\overline{\nu}(x)/\log x\to-\gamma$
as $x\to0$.
\end{enumerate}
If $F(x)$ is absolutely continuous near the origin, that is, if there
is a
measurable function~$f(x)$ such that $F(x)=\int_{0}^{x}f(u)\,\mathrm
{d}u$ for
all $x\ge0$
in a neighborhood of the origin, then
\begin{enumerate}[(S8)]
\item[(S8)] \hypertarget{a8} $\log f(x)/\log x\to\gamma-1$ as $x\to0$
\end{enumerate}
implies \textup{(S4)--(S7)}. If additionally the density $f$ is
monotone near
the origin then \textup{(S8)} is equivalent to \textup{(S4)--(S7)}.
\end{theorem}

\begin{pf}
(S3)${}\Rightarrow{}$(S1) is obvious.

(S1)${}\Leftrightarrow{}$(S2). This is clearly equivalent to
Proposition \ref{be}.

(S2)${}\Rightarrow{}$(S3).
Suppose that
$t\Phi(\mathrm{e}^{z/t})\to-\log(1-\fix(\mathrm{e}^z))$ for all
continuity points
$\mathrm{e}^z$ of~$\fix$. We know already that
$\fix(\mathrm{e}^z)=0$ for $z<0$. Moreover, if $z > 0$ then
%
\begin{equation}\label{dr}
\frac{\Phi(\mathrm{e}^{z/t})}{z/t}\to-\frac{\log(1-\fix(\mathrm
{e}^z))}{z}
\end{equation}
for all continuity points $\mathrm{e}^z$ of $\fix$. Now, the key
observation is
that the latter limit is necessarily the same for all $z>0$.
Indeed, the left-hand side of (\ref{dr}) has the form $h(z/t)$ for
some function $h$ so that if (\ref{dr})
holds for some $z>0$ then for arbitrary $z'>0$ we get, setting $t'=(z/z')t$,
\[
\lim_{t\to0}\frac{\Phi(\mathrm{e}^{z'/t})}{z'/t} =\lim_{t'\to
0}\frac{\Phi
(\mathrm{e}^{z/t'})}{z/t'} = -\frac{\log(1-\fix(\mathrm{e}^z))}{z}.
\]
Denote the limit in (\ref{dr}) by $\gamma$. As $F^*$ attains a value
in $(0,1)$, we have $\gamma\in(0,\infty)$.
Then $\fix(\mathrm{e}^z)=1-\mathrm{e}^{-\gamma z}$, that is, $F=\Pi
_\gamma$. This
completes the proof of the equivalence of (S1)--(S3).

(S4)${}\Leftrightarrow{}$(S5). This follows by setting $s=\mathrm{e}^{z/t}$
in (\ref{dr}).

(S6)${}\Rightarrow{}$(S5). For every $s\geq0$ and every $z\geq0$,
we have the decomposition
%
\begin{equation}\label{le}
\psi(s)
= s\int_0^\infty\mathrm{e}^{-sx}F(x)\,\mathrm{d}x
= \int_0^z \mathrm{e}^{-x}F(x/s)\,\mathrm{d}x+\int_z^\infty\mathrm
{e}^{-x}F(x/s)\,\mathrm{d}x.
\end{equation}
Consequently, $\lap(s)\leq F(z/s)\int_0^z \mathrm{e}^{-x}\,\mathrm
{d}x+\int_z^\infty
\mathrm{e}^{-z}\,\mathrm{d}x
=F(z/s)(1-\mathrm{e}^{-z})+\mathrm{e}^{-z}$ and
$\lap(s)\geq\mathrm{e}^{-z}\int_0^z F(x/s)\,\mathrm
{d}x+F(z/s)\mathrm{e}^{-z}$, yielding the inequalities
%
\begin{equation}\label{ol}
F(z/s)\mathrm{e}^{-z}\leq\psi(s)\leq F(z/s)(1-\mathrm
{e}^{-z})+\mathrm{e}^{-z}.
\end{equation}
Since we assume (S6), we have $F(x)=x^{\gamma+\mathrm{o}(1)}$ as
$x\to0$.
Letting $z$ be a constant on the left-hand side of \eqref{ol}, we see
that $\lap(s)\geq s^{-\gamma+\mathrm{o}(1)}$ as $s\to\infty$.
Moreover, by choosing $z=z(s)=(\log s)^2$ on the right-hand side of
\eqref{ol} we obtain
$z(s)=s^{\mathrm{o}(1)}$ and $s^\delta=\mathrm{o}(\mathrm
{e}^{z(s)})$ for any $\delta>0$, as
$s\to\infty$, so that by \eqref{ol}
\[
\psi(s)\leq F\bigl(s^{\mathrm{o}(1)-1}\bigr)+\mathrm{e}^{-(\log
s)^2} = s^{(\mathrm{o}(1)-1)(\gamma+\mathrm{o}(1))}
+ s^{-\log s}.
\]
Hence,
$\psi(s)=s^{-\gamma+\mathrm{o}(1)}$ as $s\to\infty$, which is
tantamount to (S5).

(S5)${}\Rightarrow{}$(S6). By letting $x=z/s$, it follows from
\eqref{ol} that
%
\begin{eqnarray}\label{ol2}
\frac{\mathrm{e}^z\psi(z/x)-1}{\mathrm{e}^{z}-1}\leq F(x)\leq\psi
(z/x)\mathrm{e}^{z}.
\end{eqnarray}
Now suppose that $\log\lap(s)/\log s\to-\gamma$ as $s\to\infty$.
Letting $z$ be constant on the right-hand side, we obtain that
$F(x)\leq x^{\gamma+\mathrm{o}(1)}$. Then, by choosing $z(s)=\sqrt
{\log s}$,
we see that $F(x)\geq x^{\gamma+\mathrm{o}(1)}$ and hence
$F(x)=x^{\gamma+\mathrm{o}(1)}$.

(S5)${}\Rightarrow{}$(S7). This follows from Lemma 5.17(ii) in
\cite{Kyp}.

(S7)${}\Rightarrow{}$(S5). Suppose that $\afrac{\tnu(x)}{\log x
}\to-\gamma$ as $x\to0$. Applying integration by parts, we can
write $\Phi$ as an ordinary Laplace transform:
\begin{eqnarray}\label{I}
\Phi(s)&=&\int_0^\infty\mathrm{e}^{-x}\tnu(x/s)\,\mathrm{d}x =
\int_0^K \mathrm{e}^{-x}\tnu
(x/s)\,\mathrm{d}x + \int_K^\infty
\mathrm{e}^{-x}\tnu(x/s)\,\mathrm{d}x\nonumber\\ [-8pt]\\ [-8pt]
& =&I_K(s)+J_K(s),\qquad
\mbox{say, for every } K>0.\nonumber
\end{eqnarray}
Fix an arbitrary $\varepsilon>0$. By assumption, there is an
$s_{K,\varepsilon} >K$ such that
$-\gamma-\varepsilon\le\tnu(x/s)/\log(x/s) \le-\gamma
+\varepsilon$ for all $x\in(0,K]$
and all $s\ge s_{K,\varepsilon} $. Hence,
\begin{eqnarray*}
\int_0^K \mathrm{e}^{-x}\tnu(x/s)\,\mathrm{d}x &\le&\int_0^K
\mathrm{e}^{-x}|\log(x/s)|
(\gamma+\varepsilon) \,\mathrm{d}x\\
&\le&(\gamma+\varepsilon)\biggl[\int_0^K \mathrm{e}^{-x}\log s\,
\mathrm{d}x - \int
_0^K \mathrm{e}^{-x}\log x\,\mathrm{d}x \biggr],\qquad
s>s_{K,\varepsilon}.
\end{eqnarray*}
Clearly, $\int_0^K \mathrm{e}^{-x}|\log x|\,\mathrm{d}x
<\infty$. Therefore,
\[
\limsup_{s\to\infty} I_K(s)/\log s \le(\gamma+\varepsilon)\int
_0^K \mathrm{e}^{-x}\,\mathrm{d}x
\le\gamma+\varepsilon
\]
for every $\varepsilon>0$. Thus,
%
\begin{equation}\label{II}
\limsup_{s\to\infty} I_K(s)/\log s \le\gamma.
\end{equation}
For $x>K$ we have $\tnu(x/s) \le\tnu(K/s)$ so that
$J_K(s)\le\tnu(K/s) \int_K^\infty\mathrm{e}^{-x}\,\mathrm{d}x$,
yielding
%
\begin{equation}\label{III}
\limsup_{s\to\infty} J_K(s)/\log s \le\gamma\mathrm{e}^{-K}\qquad
\mbox{for every } K>0.
\end{equation}
Letting $K\to\infty$ we obtain from (\ref{I})--(\ref{III}) that
\[
\limsup_{s\to\infty}\Phi(s)/\log s\le\gamma.
\]
The relation
$\liminf_{s\to\infty} \Phi(s)/\log s \ge\gamma$ follows along
similar lines. Altogether this proves~(S5).

(S8)${}\Rightarrow{}$(S6).
Assume that $F$ is absolutely continuous around the origin with a
density~$f$ satisfying
$\log f(x)/\log x\to\gamma-1$ as $x\to0$. Then, given an arbitrary
$\veps\in(0, \gamma)$,
we have $f(x)\leq x^{\gamma-1-\veps}$ for all sufficiently small
$x>0$. Thus,
\[
F(x)=\int_0^x f(u)\,\mathrm{d}u\leq\frac{x^{\gamma-\veps}}{\gamma
-\veps},
\]
which implies $F(x)\leq x^{\gamma-(\veps/2)}$ for sufficiently small
$x$. Similarly, it follows that $F(x)\geq x^{\gamma+(\veps/2)}$ for
small $x$, so that indeed $\lim_{x\to0}\log F(x)/\log x = \gamma$.

(S6)${}\Rightarrow{}$(S8).
Finally, suppose that $\lim_{x\to0}\log F(x)/\log x = \gamma$ and
$F$ has a
monotone density $f$ near $0$. First, let $f$ be nondecreasing at 0.
Given an arbitrary $\veps>0$, we obtain
\[
f(x)\geq\frac{1}{x}\int_0^x f(u)\,\mathrm{d}u=\frac{F(x)}{x}\geq
x^{\gamma
-1+\veps}\qquad\mbox{for small } x.
\]
Similarly,
\[
f(x)\leq\frac{\int_x^{2x}f(u)\,\mathrm{d}u}{x}\leq\frac
{F(2x)}{x}\leq
2^{\gamma-(\veps/2)}x^{\gamma-1-(\veps/2)}\qquad
\mbox{for small } x,
\]
and the right-hand side is ultimately ${\leq}x^{\gamma-1-\veps}$ as
$x\to0$. If $f$ is nonincreasing near zero
we can interchange $\leq$ and $\geq$ in the last
inequalitites.
\end{pf}

\begin{remark}
The implication ``(S6) for some $\gamma>0 \Rightarrow$
(S1)'' was already shown in \cite{SE}.
\end{remark}

\begin{remark}
Some of the above equivalences have counterparts in the theory of
regularly varying functions. In particular (S5)$\Leftrightarrow$(S6)
has the classical form of a Tauberian theorem of the type of Theorem
8.1.7 in \cite{BGT}. However, the functions there are regularly
varying while our functions are of type $x^\gamma L(x)$ with some
$L(x)=x^{\mathrm{o}(1)}$. Note that the class of regularly varying
functions is
a subclass of the class investigated here. The extra smoothness
conditions in the Karamata theory come with the reward of being able to
conclude from $f(x)=x^\gamma L(x)$, with $L$ slowly varying, that the
Laplace transform of $f$ is of the same form $s^{-\gamma} L'(1/s)$
with a precisely determined function $L'$. As opposed to this, in our
situation the exact form of the $x^{\mathrm{o}(1)}$ terms remain
unknown, but
are not needed anyway. It is also worth mentioning that for regularly
varying functions the implication (S8)$\Rightarrow$(S4)--(S7) follows
from the monotone density theorem (Theorem 8.1.8 in \cite{BGT}).
\end{remark}

\begin{remark}
Among the possible limits of $Y_t^{-t}$ as $t$ tends to zero is the somewhat
uninteresting limit $1$, which is excluded from Theorem \ref{ttt}. Loosely
speaking, this is the $\gamma=\infty$ case of the theorem. We refrain from
stating the corresponding result here.
\end{remark}

\begin{remark}
The theorem shows that the Pareto distribution is the only possible
limit distribution of $Y_t^{-t}$ as $t\to0$. This can alternatively be
deduced as follows. Note that since~$Y_t$ is a L\'evy process,
$Y_{t}\disas\sum_{k=1}^n Y_{k,t/n}$ for any $n$, where $(Y_{i,\cdot
})_{i=1,2,\ldots,k}$ are i.i.d. copies of the process $Y_\cdot$. It
follows from Proposition \ref{pro2} that if $Y_t^{-t}\dto Y^\ast$ and
$Y_{k,t}^{-t}\dto Y_k^\ast$, then (taking the limits $Y_k^\ast$ to be
independent)
\[
Y_t^{-t/n}\disas\Biggl(\sum_{k=1}^n Y_{k,t/n} \Biggr)^{-t/n}\dto
\min\{
Y_k^\ast,k=1,2,\ldots,n\}.
\]
On the other hand $Y_t^{-t/n}\dto(Y^\ast)^{1/n}$, so that $\min\{
Y_k^\ast,k=1,2,\ldots,n\}\disas(Y^\ast)^{1/n}$. Consequently,
letting ${F}^\ast(x)=\PP{Y^\ast\leq x}$ and $\overline{F}{}^\ast
(x)=1-F^\ast(x)$, we have
$\overline{F}{}^\ast(x^n)=\overline{F}{}^\ast(x)^n$ for all $n\in\sN
$. It follows that for all $q=n/m\in\sQ$ with $n,m\in\sN$,
$\overline{F}{}^\ast(x^q)=\overline F{}^\ast(x^{1/m})^n=(\overline
{F}{}^\ast(x))^q$. Hence, $F^\ast$ is a continuous function
and $\overline{F}{}^\ast(x^r)=\overline{F}{}^\ast(x)^r$ for all $r\in
\sRp$.

We next show that $F^\ast$ is strictly monotone (unless $Y^\ast=1$
a.s., which is not of interest here). We already know from Proposition
\ref{pro1} that $F^\ast$ is concentrated on $[1,\infty)$. Let
$x,y\in[1,\infty)$ with $x\not=y$ and suppose that $F^\ast
(x)=F^\ast(y)$. It then follows that $F^\ast(x^r)=F^\ast(y^r)$ for
all $r\in\sRp$, implying $F^\ast(x)=1$ constantly for $x\in
[1,\infty)$. If this is not the case, the function $g(x)=\log
\overline F{}^\ast(\mathrm{e}^x)$ is monotone decreasing on $[1,\infty)$ and
satisfies the functional equation
$g(ry)=rg(y)$ for $y\in\sRp$, identifying $g$ as $g(x)=-\gamma y$ for
some $\gamma>0$.
\end{remark}

\begin{remark}
Adding a positive drift $ct$, $c>0$ to the subordinator
$Y_t^\ast$ changes the limiting behavior dramatically, because
$Y_t^{-t}\dto Y^\ast$ implies $(ct+Y_t)^{-t}\dto1$ by
Proposition~\ref{pro2}.
\end{remark}

\begin{remark}
Suppose that $Y_t$ does not start at zero, but $Y_0\disas B$
instead, where $B$ is a~nonnegative r.v. with $q=\PP{B=0}\in(0,1)$
and $Y_t$ is of the form $Y_t=L_t+B$, where~$L_t$ is a~subordinator
independent of $B$ with $L_t^{-t} \stackrel{
d}{\rightarrow} L^*$. If $\beta$ denotes the
LST of $B$, then $\beta(u^{1/t})$ tends to $\beta(0)=1$ for
$u<1$ and to $q$ for $u>1$ as $t\to0$. Letting $\varphi$ denote the
LST of $L_t$ and $F^*$ the distribution
function of $L^*$, it follows that
\[
\lim_{t\to0}\psi_t(u^{1/t})
=\lim_{t\to0}\varphi(u^{1/t})^t\beta(u^{1/t})
=
\cases{
1-F^\ast(u), & \quad$u<1$,\cr
q(1-F^\ast(u)), &\quad$u>1$,}
\]
so that the limiting distribution of $Y_t^{-t}$ has a atom of mass
$1-q$ at 1 and an atom of mass $q$ at infinity, as expected, since by
Proposition \ref{pro1},
\[
Y_t^{-t} \stackrel{d}{\rightarrow} \min(L^*,B^*),
\]
where $L^*$ and $B^*$ are independent,
$L^*$ has a Pareto distribution and $B^*$ attains only the values 1 and
infinity.
This way we can obtain any mixture of a Pareto distribution~$\mathcal
{P}_\gamma$ and the point mass at 1 as limiting distribution (with $F^*(x)=q1_{[0,1)}(x)+(1-qx^{-1})1_{[1,\infty)}(x)$).
\end{remark}

\begin{remark}
Example 1 (in Section \ref{s2}) shows that for parametrized families
$(\psi_t)_{t>0}$ not of the infinitely divisible form $\psi(u)^t$
other interesting limit laws can occur (e.g., the shifted exponential
distribution). Thus, there may be other limit theorems and
characterizations to be explored.
\end{remark}

\section{Applications}\label{sec4}
\subsection{Explicit examples}

\begin{example}
The following distributions are all infinitely divisible (see
\cite{Sato}, Section 2.8). A close look at their distribution
functions or densities reveals that either condition (S6) or
condition (S8) can be applied so that $Y_t^{-t}$ tends in
distribution to ${\mathcal{P}}_{\gamma}$ for some $\gamma> 0$. Note
that in most cases neither explicit formulas for the convolution powers
of $F$ nor simple expressions for $\psi_t(s)$ are known.
\begin{itemize}
\item(Gamma process) The Gamma process is a standard example of a
subordinator. The density of $Y_1$ is given by
\[
f(x)=x^{{\gamma}-1}\lambda^{\gamma}\mathrm{e}^{-\lambda
x}/\Gamma({\gamma}),
\]
where $\lambda>0$ and $\gamma>0$. Obviously $f(x)\sim x^{\gamma-1}$
as $x\to0$, implying that condition (S8) holds.
\item(Weibull distribution) If $F(x)=1-\mathrm{e}^{-x^{\gamma}}$
then $F(x)\sim
x^{\gamma}$ as $x\to0$, so that in particular condition (S6) is satisfied.
\end{itemize}
For the next three distributions, the density $f(x)$ tends to some
positive constant as $x\to0$, so that condition (S8)
holds with $\gamma=1$.
\begin{itemize}
\item(Pareto-type distribution) $f(x)=\frac{a}{(1+x)^{a+1}}$, with $a>0$.
\item(F-distribution) $f(x)=\frac{\Gamma(a)\Gamma(b)}{\Gamma
(a+b)}x^{b-1}(1+x)^{-a-b}$ with $a,b>0$.
\item(Cauchy distribution on $(0,\infty)$) $f(x)=\frac{2}{\uppi
}\dfrac{1}{1+x^2}$.
\end{itemize}
\end{example}

\begin{example}[(Generalized gamma process)]
Let $\mu$ be a $\sigma$-finite
measure on $[0,\infty)$ and suppose that the L\'evy measure is given by
\[
\nu(\mathrm{d}x)=\frac{1}{x}\int_0^\infty\mathrm{e}^{-xy}\,
\mathrm{d}\mu(y)\,\mathrm{d}x.
\]
The associated L\'evy process is called a generalized gamma process
(see \cite{arthus}) and $\mu$ is the so-called \textit{Thorin measure}.
If $\mu$ is a finite measure and $\gamma=\mu([0,\infty))$, then
$\overline{\nu}(x)/\log x\rightarrow-\gamma$ as $x\rightarrow0$,
by dominated convergence. It then follows from criterion (S7) of
the theorem that $Y_t^{-t}\stackrel{d}{\rightarrow
}{\mathcal{P}}_{\gamma}$. Note that the Gamma process corresponds to
the case where $\mu$ is the Dirac measure with mass $\gamma$ at
$y=\la$.

More examples of generalized gamma processes can be found in \cite{JS}
(the complete Bernstein functions $f$ correspond to our function $\Phi
$, $\tau(\mathrm{d}s)$ corresponds to $\nu(\mathrm{d}s)/s$ and
$\rho(\mathrm{d}t)$ to $\mu
(t)\,\mathrm{d}t/t$ in our paper). For instance, if the Thorin measure
is given by
$\mu(\mathrm{d}t)=1_{(0,\gamma)}(t)\,\mathrm{d}t$ then $\Phi(s)=(x+\gamma
)\log
(x+\gamma)-x\log x-\gamma\log\gamma$ and hence $\Phi(s)/\log x\to
\gamma$. Note that indeed $\gamma=\mu([0,\infty))$. The
corresponding L\'evy measure is given by $\nu(\mathrm{d}x)=\frac
{1-\mathrm{e}^{-\gamma
x}}{x}\,\mathrm{d}x$.
\end{example}

\begin{example}[(cf. \protect\cite{SE,SL})]
Let the density of $Y_t$ be given by
$f_t(x)=\mathrm{e}^{-x}x^{-1}tI_t(x)$, where $I_t$ is the modified Bessel
function of order one. Then the Laplace exponent is given by
\[
\Phi(s)=\log\bigl(1+s-\sqrt{s^2+2s}\bigr)
\]
and since $2s(1+s-\sqrt{s^2+2s})\to1$ as $s\to\infty$ it follows that
$\Phi(s)/\log s\to1 $ and hence $Y_t^{-t}\stackrel
{d}{\rightarrow}{\mathcal{P}}_{1}$ by criterion
(S5).
\end{example}

\begin{example}
We coin the name \textit{Dickman process} for a subordinator with
L\'evy measure $\mathrm{d}\nu(x)=\gamma x^{-1}1_{(0,1]}(x)\,\mathrm
{d}x$, where $\gamma
>0$ is some parameter. The infinitely divisible distribution function
$F$ associated with $\nu$ is the generalized Dickman distribution as
defined in \cite{pwa}. This $F$ appears for example, as
\begin{itemize}
\item the distribution of a random variable $X$ satisfying $X\disas
U^{1/\gamma}(X+1)$, where $U$ is a~uniform random variable on $[0,1]$
independent of $X$,
\item the limiting distribution of $\sum_{i=1}^n (U_1U_2\cdots
U_i)^{1/\gamma}$, where $U_1,U_2,\ldots$ are independent uniform
random variables on $[0,1]$.
\end{itemize}
The name `Dickman distribution' is due to the fact that for $\gamma=1$
the density of $F$ is given by $f(x)=\mathrm{e}^{-C }\dick(x)$, where
$C$ is
Euler's constant\vadjust{\goodbreak} and $\dick$ is the generalized Dickman function. This
function is implicitly defined by $\dick(z)=1$ for $z\in[0,1]$ and
$z\dick'(z)= \dick(z-1)$ for $z>1$. Since $\tnu(x)=-\gamma\log x$
for $x$ small enough, criterion (S7) of Theorem \ref{ttt} can be
applied; we have $Y_t^{-t}\dto\pto{\gamma}$.
\end{example}

\subsection{Subordination}

\label{www} If $X_{t}$ is another subordinator with Laplace exponent
$\varphi(s)$
and $X_t$ and $Y_t$ are independent, then both subordinate processes
$A_{t}=X_{Y_{t}}$ and $B_{t}=Y_{X_{t}}$
are again subordinators. Their Laplace exponents are
\[
\phi_{A}(s)=\Phi(\varphi(s))\quad\mbox{and}\quad\phi
_{B}(s)=\varphi
(\Phi(s)),
\]
respectively. Suppose that $Y_{t}^{-t}\stackrel{d}{
\rightarrow}{\mathcal{P}}_{\gamma}$ as $t\rightarrow0$ and that
$\delta
>0 $. It follows immediately from the representations
\[
\frac{\phi_{A}(s)}{\log s}=\frac{\Phi(\varphi(s))}{\log\varphi
(s)}\frac{%
\log\varphi(s)}{\log s}\quad\mbox{and}\quad\frac{\phi
_{B}(s)}{\log s}=%
\frac{\varphi(\Phi(s))}{\Phi(s)}\frac{\Phi(s)}{\log s}
\]
and criterion (S5) that
\[
A_{t}^{-t}\stackrel{d}{\rightarrow}{\mathcal{P}}
_{\gamma\delta}\qquad\mbox{as } t\rightarrow0 \quad\iff\quad
\lim_{s\to\infty
}\frac{\log\varphi(s)}{\log{s}} = \delta
\]
and
\[
B_{t}^{-t}\stackrel{d}{\rightarrow}{\mathcal{P}}%
_{\gamma\delta}\qquad\mbox{as } t\rightarrow0 \quad\iff\quad
\lim_{s \to
\infty}\frac{\varphi(s)}{s}= \delta.
\]

\begin{example}[(Subordination with $\alpha$-stable processes)]
Suppose that
$\varphi(s)=s^\alpha$ is the Laplace exponent of an $\alpha$-stable
subordinator $X_t$. Then it follows that
\[
A_t^{-t}\stackrel{d}{\rightarrow} {\mathcal{P}}%
_{\gamma\alpha}\quad\iff\quad Y_t^{-t}\stackrel{d}{\rightarrow}{%
\mathcal{P}}_{\gamma}.
\]
For $\alpha=1$, we deduce that $B_t^{-t}\stackrel{%
d}{\rightarrow}{\mathcal{P}}_{\gamma}$ if and
only if $Y_t^{-t}%
\stackrel{d}{\rightarrow}{\mathcal{P}}_{\gamma
}$, but in
this case we just deal with the trivial case of deterministic drift
$X_t=t$ and $B_t=Y_t$.
\end{example}

\subsection{Exponential dispersion models}

For each $\theta\geq0$, we define a new L\'{e}vy measure $\nu^{%
(\theta)}$ by exponentially tilting $\nu$, that
is, we let
\[
\mathrm{d}\nu^{(\theta)}(x)=\mathrm{e}^{-\theta x}\,\mathrm{d}\nu(x).
\]
The Laplace exponent of the associated L\'{e}vy process
$Y_{t}^{(\theta)}$ is given by the difference
\[
\Phi^{(\theta)}(s)=\Phi(\theta+s)-\Phi(\theta).
\]
The new LST $\psi_{t}^{(\theta)}(s)={\mathbbl E}%
(\mathrm{e}^{-sY_{t}})$ is related to $\psi(s)$ via
\[
\psi_{t}^{(\theta)}(s)= \biggl(\frac{\psi(\theta
+s)}{%
\psi(\theta)} \biggr)^{t}.
\]
Accordingly, the distribution of $Y_{t}^{(\theta
)}$ is
given by
\[
F_{t}^{(\theta)}(\mathrm{d}x)=\frac{\mathrm{e}^{-sx} F^{t\ast
}(\mathrm{d}x)}{%
(\int_{0}^{\infty}\mathrm{e}^{-su}\,\mathrm{d}F(u))^{t}},
\]
where $F^{t\ast}$ denotes the distribution\vspace*{1pt} with LST
$(\int
_{0}^{\infty
}\mathrm{e}^{-su}\,\mathrm{d}F(u))^{t}$. The class $\{F_{t}^{(\theta
)%
},\break t\geq0,\,\theta\geq0\}$ is called an exponential dispersion model
(see \cite{SL}).

By writing
\[
\frac{\Phi^{(\theta)}(s)}{\log s}=\frac{\Phi
(\theta
+s)}{\log s}-\frac{\Phi(\theta)}{\log s}=\frac{\Phi(\theta
+s)}{\log s}%
+\mathrm{o}(1),\qquad\mbox{as } s \to\infty
\]
we see that $(Y_{t}^{(\theta)})^{-t}\rightarrow{%
\mathcal{P}}_{\gamma}$ if and only if $Y_{t}^{-t}\rightarrow
{\mathcal{P}}_{\gamma}$.

\section{A generalization}
\label{s3}

In the preceding sections, we have studied the convergence of $-t\log
Y_{t}$ to $%
X^{\ast}=\log Y^{\ast}$ as $t\rightarrow0$. In this section, we consider
the more general case
\[
tL(Y_{t})\stackrel{d}{\rightarrow}X^{\ast},\qquad
t\rightarrow0,
\]
where $L\dvtx(0,\infty)\rightarrow(-\infty,\infty)$ is some decreasing
function satisfying $\lim_{y\rightarrow0}L(y)=\infty$. Let
\[
L(\infty) \equiv
\lim_{y\rightarrow\infty}L(y)\in[-\infty,\infty)
\]
and denote by $L^{%
-1}\dvtx(-\infty,\infty)\rightarrow(0,\infty)$
the (decreasing) inverse function of $L$, with the convention that
$L^{%
-1}(x)=\infty$ for $x \le L(\infty)$, in which
case $1/L^{%
-1}(x)=0$.

The next result is the counterpart of Proposition \ref{be}, now for
the general case where $L$ is not the negative logarithm. To impose
suitable conditions on $L$, we need the definition of slow variation. The
function $L$ is called slowly varying at zero if $\lim_{x\rightarrow
0}L(\lambda x)/L(x) = 1$
for all $\lambda>0$. If this holds one can show that the
inverse function is rapidly varying at infinity, that is,
\[
\frac{L^{-1}(w)}{L^{
-1}(y)}%
\rightarrow\biggl(\frac{w}{y} \biggr)^{\infty}\equiv
\cases{
0, & \quad$w>y$, \cr
1, & \quad$w=y$, \cr
\infty,& \quad$w<y$,
}
\]
and the convergence is necessarily uniform for $w$ outside of intervals
$%
(y-\delta,y+\delta)$, $\delta>0$ (for both concepts see \cite
{BGT}). With
these prerequisites, we can show the following proposition.

\begin{proposition}
\label{P1} Suppose that $L\dvtx(0,\infty)\rightarrow(-\infty
,\infty)$ is
decreasing with $\lim_{y\rightarrow0}L(y)=\infty$ and that $L$ is slowly
varying at zero. Let $X^{\ast}$ be a random variable which is not
concentrated at one point. Then
\[
tL(Y_{t})\stackrel{d}{%
\rightarrow}X^{\ast}\qquad\mbox{as } t\rightarrow0
\]
if and only if
%
\begin{equation}\label{joco}
\psi_{t} \biggl(\frac{1}{L^{-1}(u/t)}
\biggr)\rightarrow
1-H^{\ast}(u) \qquad\mbox{as } t\rightarrow0
\end{equation}
for all continuity points $u$ of the distribution function $H^{\ast
}(u)={%
\mathbbl P}(X^{\ast}\leq u)$.
\end{proposition}

\begin{pf}
For $u<0,$ we always have $\lap_t\kla{\frac{1}{\gi\ff(u/t)}}\to1$
as $t\to0$, so we restrict ourselves to $u>0$. Let $F_t(x)=\PP
{Y_t\leq x}$ and let $H_t$ denote the distribution function of
$X_t\equiv tL(Y_t)$.
Since $H_t(x)=1-F_t(\gi\ff(x/t))$ for $t>0$, it follows that
\begin{eqnarray*}
\lap_t\biggl(\frac{1}{\gi\ff(u/t)}\biggr)
&=&\int_{0}^\infty\fexp\biggl(-\frac{y}{\gi\ff(u/t)}\biggr)\,
\mathrm{d}F_t(y)\\
&=&\int_{L(\infty)}^{\infty} \fexp\biggl(-\frac{\gi\ff
(x/t)}{\gi\ff(u/t)}\biggr)\,\mathrm{d}H_t(x).
\end{eqnarray*}
Hence, $\lap_t\kla{\frac{1}{\gi\ff(u/t)}}=\esym(\ii{t}{X_t}{u})$,
where $\ii{t}{x}{u}=\fexp\kla{-\afrac{\gi\ff(x/t)}{\gi\ff
(u/t)}}$. Since $\ff\!$ is slowly varying$\!$ at zero, it follows that $\gi
\ff$ is rapidly varying at $\infty$, implying that $\lim\!_{t\to0}\ii
{t}{x}{u}= 1_{\{x>u\}}$ for $x\not=u$. Furthermore, we obtain that
$\lim_{t\to0}\ii{t}{c_t(x)}{u}= 1_{\{x>u\}}$ for any function
$c_t(x)$ with $c_t(x)\to x$ as $t\to0$, since $x>u$ implies that
$c_t(x)>u$ eventually as $t\to0$ (and $x<u$ implies that $c_t(x)<u$
eventually).

\imp Suppose first that $X_t\dto\xix$. We can apply the continuous
mapping theorem in the form of Theorem 4.27 in \cite{kallenberg}. It
follows that
$\ii{t}{X_t}{u}\dto1_{\{x>u\}}$ for any continuity point $u$ of
$H^*$. Since $\ii{t}{x}{u}\in[0,1]$, we have $\esym(\ii
{t}{X_t}{u})\to\esym(1_{\{x>u\}})=1-\hix(u)$ by dominated convergence.

\pmi If on the other hand \eqref{joco} holds, then $\esym(\ii
{t}{X_t}{u})\to1-\hix(u)$ and
\begin{eqnarray*}
&&\bigl|\PP{X_t\leq u}-\bigl(1-\hix(u)\bigr)\bigr|\\
&&\quad\leq
|\PP{X_t\leq u}-\esym(\ii{t}{X_t}{u})|+|1-\hix(u)-\esym(\ii
{t}{X_t}{u})|\\
&&\quad=\bigl|\esym\bigl(1_{\{X_t\leq u\}}-\ii{t}{X_t}{u}\bigr
)\bigr|+|1-\hix
(u)-\esym(\ii{t}{X_t}{u})|.
\end{eqnarray*}
The second term on the right-hand side tends to zero as $t\to0$.
Regarding the first term,
for every $\veps>0$ and $\delta\in(0,u)$, we have for all $t$ large
enough (by uniform convergence for rapidly varying functions) that
\begin{eqnarray*}
\bigl|\esym\bigl(1_{\{X_t\leq u\}}-\ii{t}{X_t}{u}\bigr)\bigr|
&=&\int_0^{u-\delta} \bigl|\esym\bigl(1-\ii{t}{x}{u}\bigr)\bigr|\,\mathrm{d}H_t(x)
+\int_{u-\delta}^{u+\delta} \bigl|1_{\{x\leq u\}}-\ii{t}{x}{u}\bigr|\,\mathrm{d}H_t(x)\\
&&\quad{} +\int_{u+\delta}^\infty|\esym(\ii {t}{x}{u})|\,\mathrm{d}H_t(x)\\
&\leq&\veps\bigl(H_t(u-\delta)+1-H_t(u+\delta) \bigr)+H_t(u+\delta
)-H_t(u-\delta).
\end{eqnarray*}
Thus, if $u$ is a continuity point of $\hix$ it follows that $\abs
{\esym(1_{\{X_t\leq u\}}-\ii{t}{X_t}{u})}$
tends to zero too, yielding $\abs{\PP{X_t\leq u}-(1-\hix(u))}\to0$.
\end{pf}

We can now state the main result of this section, again for an
arbitrary subordinator~$Y_t$.
The proof follows along the lines of that of Theorem \ref{ttt}.

\begin{theorem}\label{iiie} Let $L$, $H^{\ast}$ and $X^{\ast}$ be as in Proposition
\textup{\ref{P1}} and that $\gamma>0$. Let ${\mathcal{E}}_{\gamma}$ be a~r.v. with
exponential distribution function $E_{\gamma}(x)=1-\mathrm{e}^{-\gamma x}$,
$x\geq0$. Then
the following statements are equivalent:

\begin{enumerate}
\item $t L(Y_t)\stackrel{d}{\rightarrow}
X^\ast$ as $t\to0$.

\item $t L(Y_t)\stackrel{d}{\rightarrow} {%
\mathcal{E}}_{\gamma}$ as $t\to0$.

\item $\Phi(1/s)/L(s)\to\gamma$ as $s\to0$.

\item $t\Phi(1/L^{-1}(u/t))\rightarrow\log(1-H^{\ast
}(u))$ as
$t\rightarrow0$, for all continuity points $u$ of $H^{\ast}$.
\end{enumerate}
\end{theorem}

Candidates other than $-\log x$ that satisfy the conditions of the theorem
are for example the functions $-(\log x)^{2k+1}$, $k=1,2,\ldots.$

\section*{Acknowledgement}
We are grateful to the anonymous referee for a
careful reading and valuable comments which improved the exposition of the
paper. This paper was written while Shaul Bar-Lev was a visiting professor at
the University of Osnabr\"{u}ck supported by the Mercator program of the Deutsche Forschungsgemeinschaft.

%

\printhistory


\begin{thebibliography}{12}

\bibitem{SE}
%
\begin{barticle}[mr]
\bauthor{\bsnm{Bar-Lev},~\bfnm{Shaul~K.}\binits{S.K.}} \AND
\bauthor{\bsnm{Enis},~\bfnm{Peter}\binits{P.}}
(\byear{1987}).
\btitle{Existence of moments and an asymptotic result based on a
mixture of
exponential distributions}.
\bjournal{Statist. Probab. Lett.}
\bvolume{5}
\bpages{273--277}.\break
\bid{doi={10.1016/0167-7152(87)90104-0}, issn={0167-7152}, mr={0896458}}
\end{barticle}
%
\endbibitem

\bibitem{SL}
%
\begin{barticle}[mr]
\bauthor{\bsnm{Bar-Lev},~\bfnm{Shaul~K.}\binits{S.K.}} \AND
\bauthor{\bsnm{Letac},~\bfnm{G{\'e}rard}\binits{G.}}
(\byear{2010}).
\btitle{The limiting behavior of some infinitely divisible exponential
dispersion models}.
\bjournal{Statist. Probab. Lett.}
\bvolume{80}
\bpages{1870--1874}.
\bid{doi={10.1016/j.spl.2010.08.013}, issn={0167-7152}, mr={2734253}}
\end{barticle}
%
\endbibitem

\bibitem{bertoin}
%
\begin{bbook}[mr]
\bauthor{\bsnm{Bertoin},~\bfnm{Jean}\binits{J.}}
(\byear{1996}).
\btitle{L\'evy Processes}.
\bseries{Cambridge Tracts in Mathematics}
\bvolume{121}.
\baddress{Cambridge}: \bpublisher{Cambridge Univ. Press}.
\bid{mr={1406564}}
\bptnote{check year}
\end{bbook}
%
\endbibitem

\bibitem{BGT}
%
\begin{bbook}[mr]
\bauthor{\bsnm{Bingham},~\bfnm{N.~H.}\binits{N.H.}},
\bauthor{\bsnm{Goldie},~\bfnm{C.~M.}\binits{C.M.}} \AND
\bauthor{\bsnm{Teugels},~\bfnm{J.~L.}\binits{J.L.}}
(\byear{1987}).
\btitle{Regular Variation}.
\bseries{Encyclopedia of Mathematics and Its Applications}
\bvolume{27}.
\baddress{Cambridge}: \bpublisher{Cambridge Univ. Press}.
\bid{mr={0898871}}
\end{bbook}
%
\endbibitem

\bibitem{JS}
%
\begin{barticle}[mr]
\bauthor{\bsnm{Jacob},~\bfnm{Niels}\binits{N.}} \AND
\bauthor{\bsnm{Schilling},~\bfnm{Ren{\'e}~L.}\binits{R.L.}}
(\byear{2005}).
\btitle{Function spaces as {D}irichlet spaces (about a paper by {W}.
{M}az'ya and {J}. {N}agel)}.
\bjournal{Z. Anal. Anwendungen}
\bvolume{24}
\bpages{3--28}.
\bid{doi={10.4171/ZAA/1228}, issn={0232-2064}, mr={2146549}}
\bptnote{check related}
\end{barticle}
%
\endbibitem

\bibitem{arthus}
%
\begin{barticle}[mr]
\bauthor{\bsnm{James},~\bfnm{Lancelot~F.}\binits{L.F.}},
\bauthor{\bsnm{Roynette},~\bfnm{Bernard}\binits{B.}} \AND
\bauthor{\bsnm{Yor},~\bfnm{Marc}\binits{M.}}
(\byear{2008}).
\btitle{Generalized gamma convolutions, {D}irichlet means, {T}horin measures,
with explicit examples}.
\bjournal{Probab. Surv.}
\bvolume{5}
\bpages{346--415}.
\bid{doi={10.1214/07-PS118}, issn={1549-5787}, mr={2476736}}
\end{barticle}
%
\endbibitem

\bibitem{Jorg}
%
\begin{bincollection}[auto:STB|2011/09/12|07:03:23]
\bauthor{\bsnm{J{\o}rgensen},~\bfnm{B.}\binits{B.}}
(\byear{2006}).
\btitle{Dispersion models}.
In \bbooktitle{Encyclopedia of Statistical Sciences}.
\baddress{New York}: \bpublisher{Wiley}.
\end{bincollection}
%
\endbibitem

\bibitem{kallenberg}
%
\begin{bbook}[mr]
\bauthor{\bsnm{Kallenberg},~\bfnm{Olav}\binits{O.}}
(\byear{2002}).
\btitle{Foundations of Modern Probability},
\bedition{2nd} ed.
\bseries{Probability and Its Applications (New York)}.
\baddress{New York}: \bpublisher{Springer}.
\bid{mr={1876169}}
\end{bbook}
%
\endbibitem

\bibitem{Kyp}
%
\begin{bbook}[mr]
\bauthor{\bsnm{Kyprianou},~\bfnm{Andreas~E.}\binits{A.E.}}
(\byear{2006}).
\btitle{Introductory Lectures on Fluctuations of {L}\'evy Processes with
Applications}.
\bseries{Universitext}.
\baddress{Berlin}: \bpublisher{Springer}.
\bid{mr={2250061}}
\end{bbook}
%
\endbibitem

\bibitem{pwa}
%
\begin{barticle}[mr]
\bauthor{\bsnm{Penrose},~\bfnm{Mathew~D.}\binits{M.D.}} \AND
\bauthor{\bsnm{Wade},~\bfnm{Andrew~R.}\binits{A.R.}}
(\byear{2004}).
\btitle{Random minimal directed spanning trees and {D}ickman-type
distributions}.
\bjournal{Adv. in Appl. Probab.}
\bvolume{36}
\bpages{691--714}.
\bid{doi={10.1239/aap/1093962229}, issn={0001-8678}, mr={2079909}}
\end{barticle}
%
\endbibitem

\bibitem{Sato}
%
\begin{bbook}[mr]
\bauthor{\bsnm{Sato},~\bfnm{Ken-iti}\binits{K.i.}}
(\byear{1999}).
\btitle{L\'evy Processes and Infinitely Divisible Distributions}.
\bseries{Cambridge Studies in Advanced Mathematics}
\bvolume{68}.
\baddress{Cambridge}: \bpublisher{Cambridge Univ. Press}.
\bid{mr={1739520}}
\end{bbook}
%
\endbibitem

\end{thebibliography}
\end{document}